\documentclass[11pt]{article}
\usepackage[a4paper, margin=0.7in]{geometry}
\usepackage[T1]{fontenc}

\usepackage{amsthm}
\usepackage{amsmath}
\usepackage{amssymb}
\usepackage{enumerate}
\usepackage{authblk}
\usepackage{float}

\usepackage{graphicx}
\usepackage{caption}
\usepackage{mathtools}
\usepackage[plain]{fancyref}
\usepackage{physics}
\usepackage{booktabs}
\RequirePackage{hyperref}
\usepackage{xcolor}

\restylefloat{table}

\newcommand{\R}{\mathbb{R}}  
\newcommand{\N}{\mathcal{N}}
\newcommand{\E}{\mathbb{E}}
\newcommand{\F}{\mathcal{F}}
\newcommand{\vt}{\vartheta}
\newcommand{\vtn}{\hat{\vt}_n}
\newcommand{\eps}{\varepsilon}
\newcommand{\indset}[1]{{I_{ #1 }}}
\newcommand{\ind}[1]{{\indset{\left\{ #1 \right\}}}}
\newcommand{\expec}[2][]{\mathbb{E}_{#1}\left[#2\right]}
\newcommand{\cov}[2][]{\text{COV}_{#1}\left(#2\right)}
\newcommand{\prob}[2][{}]{\mathbb{P}_{#1}\left(#2\right)}
\DeclareMathOperator*{\argmax}{arg\,max}
\newcommand{\lra}{\Leftrightarrow}
\newcommand{\del}{\partial}
\newcommand{\dprod}[1]{\left\langle #1 \right\rangle}
\newcommand{\opI}{o_{\mathbb{P}}(1)}
\newcommand{\sumin}{\sum_{i=1}^n}

\newcommand{\revised}[1]{#1}

\newtheorem{theorem}{Theorem}%
\newtheorem{corollary}[theorem]{Corollary}
\newtheorem{lemma}[theorem]{Lemma}

\theoremstyle{definition}
\newtheorem{definition}[theorem]{Definition}
\newtheorem{remark}[theorem]{Remark}

\newcommand{\fancyrefequlabelprefix}{eq}
\frefformat{plain}{\fancyrefequlabelprefix}{equation\fancyrefdefaultspacing(#1)}
\Frefformat{plain}{\fancyrefequlabelprefix}{Equation\fancyrefdefaultspacing(#1)}

\newcommand{\fancyreftheoremlabelprefix}{thm}
\frefformat{plain}{\fancyreftheoremlabelprefix}{Theorem\fancyrefdefaultspacing#1}
\Frefformat{plain}{\fancyreftheoremlabelprefix}{Theorem\fancyrefdefaultspacing#1}

\newcommand{\fancyrefcorollarylabelprefix}{cor}
\frefformat{plain}{\fancyrefcorollarylabelprefix}{Corollary\fancyrefdefaultspacing#1}
\Frefformat{plain}{\fancyrefcorollarylabelprefix}{Corollary\fancyrefdefaultspacing#1}

\newcommand{\fancyrefdefinitionlabelprefix}{def}
\frefformat{plain}{\fancyrefdefinitionlabelprefix}{Definition\fancyrefdefaultspacing#1}
\Frefformat{plain}{\fancyrefdefinitionlabelprefix}{Definition\fancyrefdefaultspacing#1}

\newcommand{\fancyreflemmalabelprefix}{lem}
\frefformat{plain}{\fancyreflemmalabelprefix}{Lemma\fancyrefdefaultspacing#1}
\Frefformat{plain}{\fancyreflemmalabelprefix}{Lemma\fancyrefdefaultspacing#1}

\hypersetup{colorlinks=true, linkcolor=blue, citecolor=blue}

\usepackage[backend=biber, style=authoryear, citestyle=authoryear, sorting=nyt, natbib, firstinits, dashed=false, maxcitenames=2]{biblatex}
\DeclareNameAlias{sortname}{family-given}
\makeatletter
\renewenvironment{table}
{\@float{table}[pb]}  %
{\end@float}
\makeatother

\title{Asymptotic properties of the MLE in distributional regression under random censoring}
\author{Gitte Kremling\thanks{Corresponding author, Email: \texttt{kremling@fh-aachen.de}}, Gerhard Dikta}
\affil{Fachhochschule Aachen, Germany}
\date{}

\begin{document}
	
	\maketitle
	
	\begin{abstract}
		Distributional regression aims to find the best candidate in a given parametric family of conditional distributions to model a given dataset. As each candidate in the distribution family can be identified by the corresponding distribution parameters, a common approach for this task is to use the maximum likelihood estimator (MLE) for the parameters. In this paper, we establish theoretical results for this estimator in case the response variable is subject to random right censoring. In particular, we provide proofs of almost sure consistency and asymptotic normality of the MLE under censoring. The \revised{empirical} behavior is illustrated by a simulation study \revised{and a real data example}.
	\end{abstract}
	
	\section{Introduction}
	\label{sec:intro}
	
	Researchers often study a response variable $Y \in \R$ together with covariates $X \in \R^p$, aiming to model the conditional distribution of $Y$ given $X$. When $Y$ represents a time-to-event outcome, observations may be right-censored: instead of $Y$, we only observe that the event time exceeds a censoring time $C$. Such censoring commonly appears in clinical studies but also arises in many other settings involving duration data.
	Under the random censoring model, each observation consists of $X$, the censored time $Z=\min(Y,C)$, and the censoring indicator $\delta=\ind{Y\le C}$, where $Y$ and $C$ are non-negative. \revised{Moreover, we assume $C$ to be conditionally independent of $Y$ given $X$. This corresponds to the standard independent censoring assumption in survival analysis and allows the censoring distribution to depend on covariates.} Given a parametric family of conditional densities with respect to some $\sigma$-finite measure $\nu$,
	\begin{equation*}
		\F = \{f: \R^p \times [0, \infty) \ni (x,y) \mapsto f(y\mid \vt, x) \in (0,\infty) \mid \vt \in \Theta\},
	\end{equation*}
	with $\Theta \subseteq \R^q$ denoting the set of admissible parameters and an i.i.d.\ sample $(X_i,Z_i,\delta_i)$, the goal is to find the parameter $\vt$ that best fits the data. In this paper, we study the asymptotic properties of the maximum likelihood estimator (MLE) $\vtn \coloneqq \argmax_{\vt \in \Theta} \ell_n(\vt; \{(x_i, z_i, d_i)\}_{i=1}^n)$ with 
	
	\begin{equation}
		\label{eq:loglik}
		\ell_n(\vt;\{(x_i, z_i, d_i)\}_{i=1}^n) \coloneqq \sumin d_i \log \left( f(z_i\mid \vt, x_i) \right) + (1-d_i) \log \left( 1 - F(z_i\mid \vt, x_i) \right),
	\end{equation}
	
	where $F(\cdot\mid \vt, \cdot)$ denotes the conditional distribution function corresponding to the density $f(\cdot\mid \vt, \cdot)$. We define $(1-d_i) \log(1-F(z_i\mid \vartheta, x_i))$ to be zero if $F(z_i\mid \vt, x_i) = 1$. 
	Note that $\ell_n$ reduces to the established log likelihood functions in \citet{stute1992strong,dikta2021bootstrap} in case the covariates are not taken into account %
	or the sample is not censored, %
	respectively. 

	\revised{The likelihood-based approach considered in this paper is widely used in survival analysis, see, for example, \citet{kalbfleisch2002statistical}. While asymptotic properties are discussed there and implementations, such as the survreg function in R, employ the same likelihood formulation, detailed proofs of consistency and asymptotic normality in the presence of covariates and censoring are not provided. The present paper fills this gap by establishing these properties under explicit regularity conditions.}
	
	There are three related problems comprising a much richer literature. 
	The first one is finding an appropriate candidate of the parametric family of conditional distributions $\F$ given an uncensored and hence fully observed sample $\{(X_i, Y_i)\}_{i=1}^n$. In a generalized linear model setting, the corresponding MLE was studied in \citet[Sec.~5.3]{dikta2021bootstrap} where its consistency and asymptotic normality were established.
	Another related problem is that of finding the best parametric fit for the marginal distribution of $Y$ under censoring but without taking any covariates into account. A proof of consistency for the MLE in this setting is given in \citet{stute1992strong}.
	Finally, we could consider the related problem of choosing an appropriate parametric model candidate for the regression function $m(x) = \E[Y\mid X=x]$ without assuming any specific conditional distribution but taking both censoring as well as the covariates into account. Corresponding parameter estimators were proposed in \citet{buckley1979linear, koul1981regression, stute1999nonlinear, heuchenne2007nonlinear}.	
	All of these approaches rely on some variations of an LSE. As we include a parametric model for the conditional distribution of the survival times in our assumptions and thereby have more information at hand, we can use the maximum likelihood method to estimate the parameters.

	The remainder of the paper is organized as follows. In \fref{sec:consis} and \ref{sec:asym_norm}, we establish the main results of this paper, that is, the consistency and asymptotic normality of the MLE $\vtn$.
	\Fref{sec:sim} presents the results of a simulation study, \revised{while a real data example is studied in \fref{sec:real}}. The appendix provides the proofs and additional simulation results.
	
	\section{Consistency}
	\label{sec:consis}
	
	To prove almost sure convergence of $\vtn$, we will use similar techniques as \citet{dikta2021bootstrap} and \citet{stute1992strong}, which study the parametric GLM without censoring and the random censorship model without covariates, respectively. Accordingly, it will be based on the work of \citet{perlman1972strong} using an extended Kullback-Leibler information.
	
	\begin{definition}[Extended Kullback-Leibler information]
		\label{def:kl_info}
		Let $\vt_1, \vt_2 \in \Theta$, \revised{$G(\cdot\mid x)$ be the conditional distribution function of the censoring variables $C$ given $X=x$} and $H$ be the distribution function of the covariates~$X$. We define
		\begin{align*}
			K_{G,H}(\vt_1, \vt_2) &\coloneqq \int \int \int_{-\infty}^c \log \left( \frac{f(y\mid \vt_1, x)}{f(y\mid \vt_2, x)} \right) f(y\mid \vt_1,x) \nu(dy)\\
			&\qquad\quad + \Big(1-F(c\mid \vt_1,x)\Big) \log \left( \frac{1-F(c\mid \vt_1, x)}{1-F(c\mid \vt_2, x)} \right) \revised{G(dc\mid x) H(dx)}
		\end{align*}
		to be the extended Kullback-Leibler information with respect to $G$ and $H$. For the sake of definiteness, we declare that $0\log(\frac{0}{0})$ equals zero. 
	\end{definition}
	
	Note that if, for every $x$, $G(\cdot\mid x)=\delta_\infty$, the Dirac measure centered at infinity, i.e.\ without any censoring, $K_{G,H}$ coincides with the extended Kullback-Leibler information without censoring
	defined in \citet{dikta2021bootstrap}. Moreover, in case $Y$ is independent of $X$, it matches the modified Kullback-Leibler information without covariates defined in \citet{stute1992strong}. %
	If both extreme cases are assumed, i.e.\ neither censoring nor covariates are involved, $K_{G,H}$ coincides with the classical Kullback-Leibler information.

	For $\vt_1, \vt_2 \in \Theta$, we define
	\begin{equation*}
		A_{\vt_1,\vt_2} \coloneqq \qty{(x,y) \,\Big|\, \frac{f(y\mid \vt_1,x)}{f(y\mid \vt_2,x)} \ne 1, \revised{G(y-\mid x)} < 1}
	\end{equation*}
	as the set of observations $(x,y)$ for which the conditional densities for $\vt_1$ and $\vt_2$ differ and which are almost surely not censored. Throughout this section, we assume the conditional density $f(y\mid\vt,x)$ to be strictly positive for \mbox{$(H \times \nu)$-a.e.} $(x,y)$ and all $\vt \in \Theta$. \revised{In fact, this is not a necessary condition. It would be enough to assume that the support of the conditional distribution of $Y$ given $X$ does not depend on $\vt$, which also allows for bounded distributions, for example.}
	The next lemma provides a fundamental result for the extended Kullback-Leibler information $K_{G,H}$.
	\pagebreak
	
	\begin{lemma}
		\label{lem:KL1}
		For $\vt_1,\vt_2 \in \Theta$, we have
		\begin{itemize}
			\item[(i)] $K_{G,H}(\vt_1, \vt_2) \in [0,\infty]$ and
			\item[(ii)] $K_{G,H}(\vt_1, \vt_2) = 0 \lra \displaystyle \int \int \indset{A_{\vt_1,\vt_2}} (x,y) f(y\mid\vt_1,x) \, \nu(dy) \, H(dx) = 0$.
		\end{itemize}
	\end{lemma}

	In the next step, we establish the relationship between the extended Kullback-Leibler distance $K_{G,H}$ and the log-likelihood function $\ell_n(\vt)$ defined in (\ref{eq:loglik}). If $\expec{\ell_1(\vt)}$ exists, it follows by the SLLN that
	\begin{align*}
		\lim_{n \to \infty} \frac{1}{n} \ell_n(\vt) &= \expec{\ell_1(\vt)}\\
		&= \int \expec{\ind{Y \le C} \log(f(Y\mid \vt, x)) + \ind{Y > C} \log(1-F(C\mid \vt, x))\mid X=x} H(dx)\\
		&= \int \int \bigg[ \int_{-\infty}^c \log(f(y\mid \vt, x)) f(y\mid \vt_0,x) \nu(dy)\\
		&\qquad\qquad\quad+ (1-F(c\mid \vt_0,x)) \log(1-F(c\mid \vt, x)) \bigg] \revised{G(dc\mid x)} H(dx)\\
		&\eqqcolon L_{G,H}(\vt_0, \vt).
	\end{align*}
	Further, if $|L_{G,H}(\vt_1, \vt_1)| < \infty$, we have 
	$K_{G,H}(\vt_1, \vt_2) = L_{G,H}(\vt_1, \vt_1) - L_{G,H}(\vt_1, \vt_2).$
	The next lemma gives conditions for $L_{G,H}(\vt_0, \,\cdot\,)$ to attain its unique maximum at $\vt_0$. In the following theorem, we will use several convergence theorems to transfer the result to $\ell_n$.
	
	\begin{lemma}
		\label{lem:KL2}
		If
		\begin{itemize}
			\item[(i)] $|L_{G,H}(\vt_0, \vt_0)| < \infty$ and
			\item[(ii)] $\prob{A_{\vt_0,\vt}} > 0$ for all $\vt \in \Theta \backslash \{ \vt_0 \}$,
		\end{itemize}
		then the function $L_{G,H}(\vt_0, \cdot)$ has a unique maximum at $\vt_0$.
	\end{lemma}
	
	Now, we have all the tools at hand to provide a consistency result for the MLE.
	
	\begin{theorem}
		\label{thm:consis}
		Under the assumptions of \fref{lem:KL2} and if the density functions are continuous in $\vt$ for all $x \in \R$, $\Theta$ is compact and
		\begin{enumerate}[(i)]
			\item
			\label{cond:sup_l1_int}
			for all $\vt^* \in \Theta$ there exists an open neighborhood $V^* = V(\vt^*)$ of $\vt^*$ such that 
			$ \expec{\sup_{\vt \in V^*} \ell_1(\vt)} < \infty, $
		\end{enumerate}
		then $\hat{\vt}_n \to \vt_0$ as $n \to \infty$, wp1, where $\hat{\vt}_n$ is the maximum of $\ell_n(\,\cdot\,)$.
	\end{theorem}
	\noindent
	Note that the conditions of the preceding theorem are equivalent to the ones in \citet[Theorem 5.52]{dikta2021bootstrap}, establishing the strong consistency of the MLE in case of a parametric GLM without any censoring being involved. \revised{Similarly to \cite[Corollary 5.53]{dikta2021bootstrap}, the assumption of the parameter space $\Theta$ being compact can be replaced by the weaker condition: there exists a compact set $C \subseteq \Theta$ such that $\vt_0 \in C$ and $\expec{\sup_{\vt \in \Theta\backslash C} \ell_1(\vt) - \ell_1(\vt_0)} < 0$.}
	
	\section{Asymptotic normality}
	\label{sec:asym_norm}
	
	Having established the consistency of the MLE, we now want to study its asymptotic distribution.
	For any distribution function $F$, let $\bar{F}(x) = 1 - F(x)$ represent the corresponding tail distribution function. Further, let $D_r \equiv \frac{\del}{\del \vt_r}$ and $D_{r,s} \equiv \frac{\del^2}{\del \vt_r \del \vt_s}$ denote the differential operators of first and second order with respect to the $r$-th or $r$-th and $s$-th component of $\vt$, respectively, and $D(\cdot) = (D_1(\cdot), \dots, D_q(\cdot))^T$. The following theorem establishes the asymptotic normality of the MLE.
	
	\begin{theorem}
		\label{thm:asym_norm}
		If
		\begin{enumerate}[(i)]
			\item \label{cond:mle_norm__reg1}
			$\ell_1(\vt)$ has continuous second derivatives with respect to $\vt$ and there exists an open neighborhood $V \subseteq \Theta$ of $\vt_0$ such that 
			$ \expec{\sup_{\vt \in V} |D_{r,s} \ell_1(\vt\mid X,Z,\delta)|} < \infty $
			for $H$-a.e.\ $x$ and for all $1 \le r,s \le q$,
			
			\item \label{cond:mle_norm__reg2}
			$ \int \sup_{\vt \in V} \left| D_r f(y\mid \vt,x) \right| \nu(dy) < \infty $
			for all $1 \le r \le q$ and $V$ from (\ref{cond:mle_norm__reg1}),

			\item \label{cond:mle_norm__reg3}
			for all $1 \le r,s \le q$, it holds
			$$\int \int \Bigg[ \bigg[ \int_{-\infty}^c D_{r,s}f(y\mid \vt_0,x) \nu(dy) \bigg] - D_{r,s}F(c\mid \vt_0,x) \Bigg] \revised{G(dc\mid x) H(dx)} = 0,$$ 
			
			\item \label{cond:mle_norm__sigma}
			the matrix
			\begin{align*}
				\Sigma &\coloneqq \expec{\revised{\bar{G}(Y-\mid X)} D(\tilde{\ell}_1(\vt_0)) \big(D(\tilde{\ell}_1(\vt_0)))^T} \\
				&\qquad + \expec{\big(\bar{F}(C\mid \vt_0,X)\big)^{-1} \expec{\ind{Y \le C} D(\tilde{\ell}_1(\vt_0))\mid X,C} \expec{\ind{Y \le C} \big(D(\tilde{\ell}_1(\vt_0)))^T\mid X,C}} 
			\end{align*}
			exists and is positive definite, where
			$ \tilde{\ell}_n(\vt_0) = \sum_{i=1}^n \log(f(Y_i\mid \vt_0,X_i))$
			denotes the log likelihood function without any censoring being involved,	
			
			\item \label{cond:mle_norm__consis}
			$\vtn$ converges in probability to $\vt_0$,
		\end{enumerate}
		then $\sqrt{n}(\vtn-\vt_0) \Rightarrow \mathcal{N}_q(0, \Sigma^{-1})$.
	\end{theorem}
	\noindent
	Again, the conditions of the theorem are similar to the corresponding result for parametric GLMs without censoring as stated in \citet[Theorem 5.55]{dikta2021bootstrap}.
	
	\begin{remark}
		If $F(C\mid \vt_0, X) = 1$, the value of the covariance matrix $\Sigma$ is undefined, so this case has to be handled separately. In the proof of \fref{thm:asym_norm}, the second summand of $\Sigma$ is  originally derived as
		\begin{equation*}
			\Sigma_2 = \expec{\frac{(1-\delta)^2}{\qty(1-F(C\mid \vt_0,X))^2} \expec{\ind{Y \le C} D(\tilde{\ell}_1(\vt_0))\mid X,C} \Big(\expec{\ind{Y \le C} D(\tilde{\ell}_1(\vt_0))\mid X,C}\Big)^T}.
		\end{equation*}
		Since $\delta = 1$ wp1 if $F(C\mid \vt_0, X) = 1$, we will define $\Sigma_2$ to be zero in this case.
	\end{remark}
	
	\begin{remark}
		Condition (\ref{cond:mle_norm__reg2}) is needed to guarantee that the order of integration and differentiation can be interchanged. %
		A similar assumption for the second derivative, i.e.
		$$ \int_{-\infty}^c D_{r,s} f(y\mid \vt_0,x) \nu(dy) = D_{r,s} \int_{-\infty}^c f(y\mid \vt_0,x) \nu(dy) $$
		for all $1 \le r,s \le q$, $H$-a.e.\ $x$ and $G$-a.e.\ $c$, is sufficient for condition (\ref{cond:mle_norm__reg3}) to hold.
	\end{remark}
	
	The following corollary gives sufficient conditions for assumption (\ref{cond:mle_norm__sigma}) in \fref{thm:asym_norm} to hold.%
	
	\begin{corollary}
		\label{cor:asym_norm}
		If assumptions (\ref{cond:mle_norm__reg1})-(\ref{cond:mle_norm__reg3}) and (\ref{cond:mle_norm__consis}) of the previous theorem hold and
		\begin{enumerate}[(i)]
			\setcounter{enumi}{5}
			\item \label{cond:mle_norm__censoring}
			$\revised{G(Y-\mid X)} < 1$ wp1
			\item \label{cond:mle_norm__sigmatilde}
			the Fisher-information matrix of the uncensored case
			$ \tilde{\Sigma} = \mathbb{E}[D(\tilde{\ell}_1(\vt_0))(D(\tilde{\ell}_1(\vt_0)))^T] $
			exists and is positive definite,
		\end{enumerate}
		then $\sqrt{n}(\vtn-\vt_0) \Rightarrow \mathcal{N}_q(0, \Sigma^{-1})$.
	\end{corollary}
	\noindent
	Assumption \eqref{cond:mle_norm__censoring} in the corollary above 
	ensures that there are no values of $Y$ that can never be observed because they get censored wp1. Note that this is a sufficient, not necessary, condition for the asymptotic normality of the MLE.
	
	\revised{
		\begin{remark}
			\label{rem:est_sig}
			In practice, the covariance matrix $\Sigma = \cov{D(\ell_1(\vt_0))} = -\expec{D^2(\ell_1(\vt_0))}$, see \eqref{eq:cov_exp}, can be estimated via a plug-in approach, %
			yielding $\hat{\Sigma} \coloneqq \frac{1}{n} D^2 \ell_n(\vtn)$. 
		\end{remark}
	}
	
	The following corollary provides an asymptotically linear representation for the %
	MLE $\vtn$.
	\pagebreak
	
	\begin{corollary}
		\label{cor:lin_rep}
		With 
		$\xi(X, Z, \delta, \vt_0) \coloneqq \Sigma^{-1} D\big(\ell_1(\vt_0\mid X,Z,\delta)\big)$
		and under the assumptions of \fref{thm:asym_norm}, it follows that
		\begin{enumerate}[(i)]
			\item $\sqrt{n}(\vtn - \vt_0) = n^{-1/2} \sum_{i=1}^n \xi(X_i, Z_i, \delta_i, \vt_0) + o_\mathbb{P}(1)$,
			\item $\expec{\xi(X, Z, \delta, \vt_0)} = 0$,
			\item $\expec{\xi(X, Z, \delta, \vt_0) \xi^T(X, Z, \delta, \vt_0)}$ exists and is positive definite.
		\end{enumerate}
	\end{corollary}
	
	\revised{The theoretical analysis presented in this paper relies on the conditional independence of $C$ and $Y$ given $X$. If censoring were informative, so that the distribution of $C$ depends on $Y$ even conditional on $X$, an approach similar to that proposed by \cite{dikta1998semiparametric} could be employed. The main idea is to introduce another parametric model $m(\theta, x, z)$ for $\expec{\delta \mid X=x, Z=z}$ and replace the censoring indicators $d_i$ in the likelihood function \eqref{eq:loglik} by $m(\theta, X_i, Z_i)$. %
		Extending the present results to such settings is beyond the scope of this paper.}\\[-\baselineskip]
	
	\section{Simulation study}
	\label{sec:sim}
	
\begin{figure}[h]
	\includegraphics[width=\textwidth]{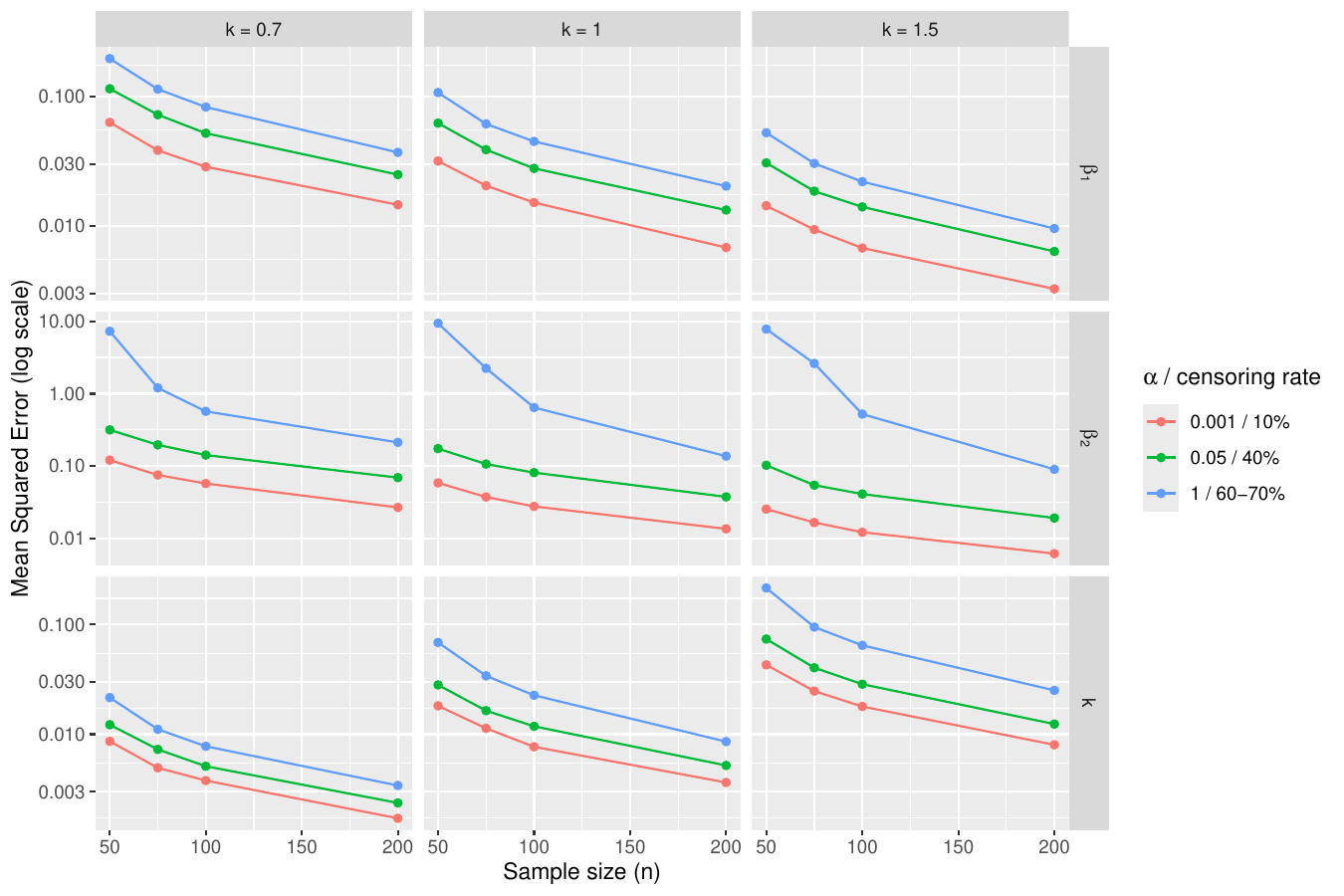}
	\caption{\revised{Accuracy of the MSE in our simulation study using different values for $\alpha$ (controlling the censoring rate) and $k$ (determining the shape of the Weibull distribution).}}
	\label{fig:sim}
	\end{figure}

	In this section, we present the results of a simulation study illustrating the asymptotic behavior of the MLE. \revised{For that, we use a two-dimensional covariate vector with $X_1 \sim \N(0,1)$ and $X_2 \sim \text{Bernoulli}(0.5)$. For the conditional distribution of $Y$ given $X$, we consider a Weibull distribution with shape parameter $k$ and scale chosen such that $\expec{Y|X=x} = \exp(\beta^Tx)$. The censoring times are sampled from an exponential distribution with $\expec{C|X=x} = \alpha \exp(\gamma^T x)$. We set $\beta = (1,2)$ and $\gamma = (2,2)$, and vary $\alpha \in \{0.001, 0.05, 1\}$ and $k \in \{0.7, 1, 1.5\}$ to study how the quality of the estimates is affected by different censoring rates and distribution shapes. Note that the choice of $\alpha$ directly controls the censoring rate. The distribution of censoring rates under different values of $\alpha$ is illustrated in Figure~\ref{fig:cens_rates} in Appendix~\ref{sec:app_sim}.} According to these distributions, we generate an i.i.d.\ sample of observations $\{(X_i, Z_i, \delta_i)\}_{i=1}^n$ with $n = \revised{50, 75,} 100, 200, 500$ and calculate the MLE for $\vt = (\beta, k) \in \R^3$. For each choice of $\alpha$ and $k$, the procedure is \revised{repeated $500$ times} to determine the %
	mean squared error (MSE) of the estimators. The results are presented in Figure \ref{fig:sim}.%
	
	\revised{In all nine combinations of $\alpha$ and $k$}, the MSE of all three parameter estimates decreases substantially with larger sample sizes. \revised{The $y$-axis in the plots is logarithmically scaled}. As expected, higher censoring rates \revised{(corresponding to larger values of $\alpha$)} lead to higher MSEs.
	\revised{Nevertheless, even in the most extreme scenario with a censoring rate around $60$--$70\%$, the estimates remain reasonably accurate for $n=500$. We also examined the performance of the estimator for very small sample sizes and found that it remains unbiased. However, as anticipated, high censoring rates in this extreme case result in substantially increased variance. A detailed illustration of the results for $n=25$ is provided in Figure~\ref{fig:n25} in Appendix~\ref{sec:app_sim}.}
	\vspace*{\baselineskip}
	
	\section{Real data example}
	\label{sec:real}
	
	\revised{In this section, we demonstrate the application of the MLE to a real data example. Specifically, we consider the Mayo Clinic Primary Biliary Cholangitis (PBC) Data that is available via the survival package in R and was also studied in \cite{therneau2000modeling}. Primary biliary cholangitis is an autoimmune disease leading to destruction of the small bile ducts in the liver. The dataset contains observations from 418 PBC patients who participated in a trial conducted between 1974 and 1984 at the Mayo Clinic. The censored times $Z_i$ correspond to the time in days between registration and the earlier of death, transplantation, or study analysis. We consider death as the event time $Y_i$ and all other values as censored. In the vector of covariates $X_i \in \R^4$, we include the patient's age (in years) and sex (0 for masculine and 1 for feminine), as well as the logarithms of serum bilirubin (in mg/dl) and serum albumin (in g/dl). As parametric family for the conditional distribution of $Y$ given $X$, we consider a Weibull distribution with shape parameter $k$ and scale chosen to satisfy $\expec{Y|X=x} = \exp(\beta_0 + \beta^T X)$. We use the MLE to estimate the parameter vector $\vt = (\beta_0, \beta, k) \in \R^6$. 
	
	The results are summarized in Table~\ref{tab:PBC}. In our implementation, the likelihood function is numerically maximized using the optim routine in R, and the confidence intervals are computed based on the plug-in estimator $\hat{\Sigma}$ defined in Remark~\ref{rem:est_sig}. The resulting estimates and confidence intervals closely agree with those obtained using the standard survreg function in R. To assess the plausibility of these confidence intervals, Table~\ref{tab:PBC} also reports the mean and corresponding quantiles from a non-parametric bootstrap based on 500 resamples of the observed data, which generally match the asymptotic intervals well. For $\beta_2$, the bootstrap interval is slightly wider. Analogous results for $\beta_0$ and the shape parameter $k$ are provided in Table~\ref{tab:PBC2} in Appendix~\ref{sec:app_pbc}.}
	
	\begin{table}
	\centering
	\caption{\revised{Estimates and 95\% confidence intervals for the linear coefficients $\beta$, corresponding to the patient's (1) age, (2) sex, (3) logarithm of serum bilirubin and (4) logarithm of serum albumin.}}
	\label{tab:PBC}
	\centering
	\revised{
		\begin{tabular}[t]{lcccc}
			\toprule
			Method & $\beta_1$ & $\beta_2$ & $\beta_3$ & $\beta_4$\\
			\midrule
			Ours & -0.03 (-0.04, -0.02) & 0.06 (-0.25, 0.38) & -0.64 (-0.75, -0.52) & 2.14 (1.34, 2.94)\\
			Survreg & -0.03 (-0.04, -0.02) & 0.06 (-0.25, 0.38) & -0.64 (-0.75, -0.52) & 2.14 (1.34, 2.93)\\
			Bootstrap & -0.03 (-0.04, -0.02) & 0.07 (-0.34, 0.45) & -0.64 (-0.77, -0.53) & 2.15 (1.35, 2.98)\\
			\bottomrule
		\end{tabular}
	}
	\end{table}
	
	\pagebreak
	\printbibliography
	
	\appendix
	\section*{Appendix}
	\section{Proofs of main theorems}
	
	\begin{proof}[Proof of \fref{thm:consis}]
		First, note that the continuity of the density functions and compactness of $\Theta$ ensure the existence of $\hat{\vt}_n$. Let $V$ be an arbitrary neighborhood of $\vt_0$ and $U = \Theta \backslash V$. As stated in \citet{perlman1972strong}, it is sufficient to prove that
		\begin{equation}
			\prob{\limsup_{n \to \infty} \sup_{\vt \in U} \frac{1}{n} \Big(\ell_n(\vt) - \ell(\vt_0)\Big) < 0} = 1.
		\end{equation}
		Note that we might have measurability issues because $U$ is uncountable. In this case, we would use the inner probability measure.
		
		Let $\vt^* \in U$ and $\{ V_\eps (\vt^*) \}$ be open neighborhoods of $\vt^*$ with $\cap_{\eps>0} V_\eps(\vt^*) = \{ \vt^* \}$. Choose $\eps > 0$ such that $V_\eps(\vt^*) \subset V^*(\vt^*)$ and $0 < M < \infty$. Then we have
		\begin{align*}
			\sup_{\vt \in V_\eps(\vt^*)} \frac{1}{n} \Big(\ell_n(\vt) - \ell_n(\vt_0)\Big) &\le \frac{1}{n} \sum_{i=1}^n \sup_{\vt \in V_\eps(\vt^*)} \max \{\delta_i \log{f(Z_i\mid \vt,X_i)} + (1-\delta_i) \log{(1 - F(Z_i\mid \vt, X_i))}, -M \}\\
			&\quad - \frac{1}{n} \sum_{i=1}^n \delta_i \log{f(Z_i\mid \vt_0, X_i)} + (1-\delta_i) \log{(1 - F(Z_i\mid \vt_0, X_i))}.
		\end{align*}
		Applying a series of convergence theorems will establish that this last expression is negative wp1.
		
		First, we apply the SLLN, 
		giving us 
		\begin{align*}
			&\limsup_{n \to \infty} \sup_{\vt \in V_\eps(\vt^*)} \frac{1}{n} \Big(\ell_n(\vt) - \ell_n(\vt_0)\Big)\\
			&\quad \le \expec{\sup_{\vt \in V_\eps(\vt^*)} \max \{\delta \log{f(Z\mid \vt, X)} + (1-\delta) \log{(1 - F(Z\mid \vt, X))}, -M \}} - L_{G,H}(\vt_0, \vt_0). 
		\end{align*}
		Note that assumption (\ref{cond:sup_l1_int}) guarantees that the expectations exist.
		In a next step, we apply Lebesgue's Dominated Convergence Theorem (DCT). Again, assumption~(\ref{cond:sup_l1_int}) ensures that there is an integrable dominating function. 
		So we have
		\begin{align*}
			&\lim_{\eps \to 0} \expec{\sup_{\vt \in V_\eps(\vt^*)} \max \{\delta \log{f(Z\mid \vt, X)} + (1-\delta) \log{(1 - F(Z\mid \vt, X))}, -M \}} - L_{G,H}(\vt_0, \vt_0)\\
			&\quad = \expec{\max \{\delta \log{f(Z\mid \vt^*, X)} + (1-\delta) \log{(1 - F(Z\mid \vt^*, X))}, -M \}} - L_{G,H}(\vt_0, \vt_0).
		\end{align*}
		Finally, we apply the Fatou-Lebesgue Theorem, which is applicable since $\left( \ell_1(\vt^*) \right)^+$ is an integrable function dominating the integrand. %
		It follows that
		\begin{align*}
			&\limsup_{M \to \infty} \expec{\max \{\delta \log{f(Z\mid \vt^*, X)} + (1-\delta) \log{(1 - F(Z\mid \vt^*, X))}, -M \}} - L_{G,H}(\vt_0, \vt_0)\\
			&\quad \le L_{G,H}(\vt_0, \vt^*) - L_{G,H}(\vt_0, \vt_0) < 0,
		\end{align*}
		where the last inequality is a direct consequence of \fref{lem:KL2}. Altogether, this establishes that for small enough $\eps$
		\begin{align*}
			\limsup_{n \to \infty} \sup_{\vt \in V_\eps(\vt^*)} \frac{1}{n} \Big(\ell_n(\vt) - \ell_n(\vt_0)\Big) < 0.
		\end{align*}
		Since $U$ is compact, there exist $\eps_i$ and $\vt^*_i$ such that $U \subset \cup_{i=1}^{m} V_{\eps_i}(\vt^*_i)$ which completes the proof.
	\end{proof}

	\begin{proof}[Proof of \fref{thm:asym_norm}]
		For the first part of the proof, we closely follow \citet{dikta2021bootstrap}. Define $s_n(\vt) \coloneqq D(\ell_n(\vt))$ and note that
		\begin{equation*}
			s_n(\vtn) - s_n(\vt_0) = \qty(\int_0^1 D\big(s_n(\vt_0 + t(\vtn - \vt_0))\big) dt) (\vtn-\vt_0).
		\end{equation*}
		Now, we substitute the integral by $D(s_n(\vt_0))$. For that, define
		\begin{equation*}
			\Delta_n \coloneqq \int_0^1 D \big(s_n(\vt_0 + t(\vtn - \vt_0))\big) dt - D(s_n(\vt_0))
		\end{equation*}
		and $B_\eps \coloneqq \{\vt : ||\vt-\vt_0|| \le \eps \}$. 
		Using Markov's inequality, we then have for the $r$-th and $s$-th component of $\Delta_n$, denoted by $\Delta_n^{(r,s)}$, that
		\begin{align*}
			\prob{|\Delta_n^{(r,s)}/n| > \tilde{\eps}} &\le \prob{\frac{1}{n} \int_0^1 \big|D_{r,s} \ell_n(\vt_0+t(\vtn-\vt_0)) - D_{r,s} \ell_n(\vt_0) \big| dt > \tilde{\eps}, \vtn \in B_\eps} + \prob{\vtn \notin B_\eps}\\
			&\le \prob{\frac{1}{n} \sup_{\vt \in B_\eps} \big|D_{r,s} \ell_n(\vt) - D_{r,s} \ell_n(\vt_0) \big| > \tilde{\eps}} + \prob{\vtn \notin B_\eps}\\
			&\le \tilde{\eps}^{-1} \expec{\frac{1}{n} \sup_{\vt \in B_\eps} \big|D_{r,s} \ell_n(\vt) - D_{r,s} \ell_n(\vt_0)\big|} + \prob{\vtn \notin B_\eps}\\
			&\le \tilde{\eps}^{-1} \expec{\sup_{\vt \in B_\eps} \big| D_{r,s} \ell_1(\vt) - D_{r,s} \ell_1(\vt_0)\big|} + \prob{\vtn \notin B_\eps}.
		\end{align*}
		Since, by (\ref{cond:mle_norm__consis}), $\vtn$ converges in probability to $\vt_0$, the second term on the right-hand side converges to zero. By the continuity assumption of (\ref{cond:mle_norm__reg1})
		$
		\lim_{\eps \to 0} \sup_{\vt \in B_\eps} \big|D_{r,s} \ell_1(\vt) - D_{r,s} \ell_1(\vt)\big| = 0$
		and hence, by assumption~(\ref{cond:mle_norm__reg1}) and the DCT, the first term on the right-hand side can be made arbitrarily small as well. Altogether, we obtain $n^{-1} \Delta_n = o_{\mathbb{P}_n}(1)$ and since $s_n(\vtn) = 0$, the initial equality becomes
		\begin{align}
			\label{eq:lemma_rep}
			-n^{-1/2} s_n(\vt_0) = (n^{-1} D(s_n(\vt_0)) + o_{\mathbb{P}}(1))n^{1/2}(\vtn-\vt_0).
		\end{align}
		
		The next steps are to apply the CLT to $s_n(\vt_0)$, the SLLN to $D(s_n(\vt_0))$ and finally \citet[Lemma 5.54]{dikta2021bootstrap}. For that, we first need to determine the expected values of the two terms and the variance of the former. We have
		\begin{align*}
			s_n(\vt_0) &= D(\ell_n(\vt_0))\\
			&= \sum_{i=1}^n \delta_i D(\log(f(Z_i\mid \vt_0, X_i))) + \frac{1-\delta_i}{1-F(C_i\mid \vt_0,X_i)} D(1-F(C_i\mid \vt_0,X_i))\\
			&= \sum_{i=1}^n \delta_i D(\log(f(Z_i\mid \vt_0, X_i))) - \frac{1-\delta_i}{1-F(C_i\mid \vt_0,X_i)} \int_{-\infty}^{C_i} D(f(y\mid \vt_0,X_i)) \nu(dy)\\
			&= \sum_{i=1}^n \delta_i D(\log(f(Z_i\mid \vt_0, X_i))) \\
			&\qquad- \frac{1-\delta_i}{1-F(C_i\mid \vt_0,X_i)} \int \ind{y \le C_i} D(\log(f(y\mid \vt_0,X_i))) f(y\mid \vt_0,X_i) \nu(dy).
		\end{align*}
		Note that we were allowed to change the order of differentiation and integration by DCT and assumption~(\ref{cond:mle_norm__reg2}).
		In the subsequent steps, we will repeatedly use the identity
		\begin{align*}
			\int \varphi(X,y,C) f(y\mid \vt,X) \nu(dy) = \expec{\varphi(X,Y,C)\mid X, C},
		\end{align*}
		which can be easily verified using \revised{the conditional independence of $C$ and $Y$ given $X$}. In particular, this yields
		\begin{align*}
			1-F(C\mid \vt_0,X) = \int \ind{y>C} f(y\mid \vt_0,X) \nu(dy) = \expec{\ind{Y>C}\mid X,C}.
		\end{align*}
		With these representations, we get
		\begin{align*}
			\expec{s_1(\vt_0)} %
			&= \expec{\ind{Y \le C} D(\log(f(Z\mid \vt_0, X)))} - \expec{\frac{\ind{Y>C}}{\expec{\ind{Y>C}\mid X,C}} \expec{\ind{Y \le C} D(\log(f(Y\mid \vt_0,X))) \mid  X, C}}\\
			&= \expec{\ind{Y \le C} D(\tilde{\ell}_1(\vt_0))} - \expec{\frac{\expec{\ind{Y>C}\mid X,C}}{\expec{\ind{Y>C}\mid X,C}} \expec{\ind{Y \le C} D(\tilde{\ell}_1(\vt_0)) \mid  X, C}}\\
			&= 0
		\end{align*}
		and hence also $\expec{s_n(\vt_0)} = 0$.
		
		Next, we consider $\cov{s_n(\vt_0)}$ and $\expec{D(s_n(\vt_0))}$. For readability purposes, we suppress the appearance of parameters in $f$ and $F$. Note that for any $1 \le r,s \le q$
		\begin{equation*}
			D_r \big(\ell_1(\vt_0)\big) = D_r \big(\delta \log(f) + (1-\delta)\log(1-F)\big) = \frac{\delta}{f} D_r(f) - \frac{1-\delta}{1-F} D_r(F)
		\end{equation*}
		and, since $\delta^2 = \delta$ and $\delta(1-\delta)=0$,
		\begin{align*}
			D_{r,s}(\ell_1(\vt)) %
			&= -\frac{\delta^2}{f^2} D_r(f) D_s(f) - \frac{(1-\delta)^2}{(1-F)^2} D_r(F) D_s(F) + \frac{\delta}{f} D_{r,s}(f) - \frac{1-\delta}{1-F} D_{r,s}(F)\\
			&= \qty(\frac{\delta}{f} D_r(f) - \frac{1-\delta}{1-F} D_r(F)) \qty(-\frac{\delta}{f} D_s(f) + \frac{1-\delta}{1-F} D_s(F)) + \frac{\delta}{f} D_{r,s}(f) - \frac{1-\delta}{1-F} D_{r,s}(F)\\
			&= -D_r(\ell_1(\vt)) D_s(\ell_1(\vt)) + \frac{\delta}{f} D_{r,s}(f) - \frac{1-\delta}{1-F} D_{r,s}(F).
		\end{align*}
		Thus, the expected value is given by
		\begin{align*}
			&\expec{D_{r,s} (\ell_1(\vt_0))}\\
			&\quad= -\expec{D_r(\ell_1(\vt_0)) D_s(\ell_1(\vt_0))} + \expec{\frac{\ind{Y \le C}}{f(Y\mid \vt_0,X)} D_{r,s}(f(Y\mid \vt_0,X)) - \frac{\ind{Y> C}}{1-F(C\mid \vt_0,X)} D_{r,s}(F(C\mid \vt_0,X))} \\
			&\quad= -\cov{D_r(\ell_1(\vt_0)), D_s(\ell_1(\vt_0))}  \\
			&\qquad+ \expec{\frac{\ind{Y \le C}}{f(Y\mid \vt_0,X)} D_{r,s}(f(Y\mid \vt_0,X))} - \expec{\frac{\expec{\ind{Y> C}\mid X,C}}{\expec{\ind{Y> C}\mid X,C}} D_{r,s}(F(C\mid \vt_0,X))} \\
			&\quad= -\cov{D_r(\ell_1(\vt_0)), D_s(\ell_1(\vt_0))} + \int \int \int_{-\infty}^c D_{r,s}(f(y\mid \vt_0,x)) \nu(dy) - D_{r,s}(F(c\mid \vt_0,x)) \revised{G(dc\mid x) H(dx)},\\
			&\quad= -\cov{D_r(\ell_1(\vt_0)), D_s(\ell_1(\vt_0))},
		\end{align*}
		where assumption (\ref{cond:mle_norm__reg3}) guarantees the last equality. In other words, we have
		\begin{align}
			\label{eq:cov_exp}
			\expec{D(s_n(\vt_0))} = -\cov{s_n(\vt_0)}.
		\end{align}
		For the covariance matrix, we get
		\begin{align*}
			&\cov{s_1(\vt_0)} = \expec{s_1(\vt_0) \big(s_1(\vt_0)\big)^T}\\			
			&\quad= \expec{\delta^2 D(\log(f(Z\mid \vt_0,X))) \big(D(\log(f(Z\mid \vt_0,X)))\big)^T}\\
			&\qquad\quad + \mathbb{E}\left[\frac{(1-\delta)^2}{\qty(1-F(C\mid \vt_0,X))^2} \Big(\int \ind{y \le C} D(\log(f(y\mid \vt_0,X))) f(y\mid \vt_0,X) \nu(dy)\Big) \right.\\
			&\qquad\qquad\qquad\qquad\qquad\qquad\quad~ \left. \cdot \Big(\int \ind{y \le C} D(\log(f(y\mid \vt_0,X))) f(y\mid \vt_0,X) \nu(dy)\Big)^T\right]\\
			&\quad= \expec{\ind{Y \le C} D(\tilde{\ell}_1(\vt_0)) \big(D(\tilde{\ell}_1(\vt_0))\big)^T}\\
			&\qquad\quad + \expec{\frac{\ind{Y > C}}{\qty(\expec{\ind{Y > C}\mid X,C})^2} \Big(\expec{\ind{Y \le C} D(\tilde{\ell}_1(\vt_0))\mid X,C}\Big) \cdot \Big(\expec{\ind{Y \le C} D(\tilde{\ell}_1(\vt_0))\mid X,C}\Big)^T}\\
			&\quad= \expec{\ind{Y \le C} D(\tilde{\ell}_1(\vt_0)) \big(D(\tilde{\ell}_1(\vt_0))\big)^T}\\
			&\qquad\quad + \expec{\frac{1}{\expec{\ind{Y > C}\mid X,C}} \Big(\expec{\ind{Y \le C} D(\tilde{\ell}_1(\vt_0))\mid X,C}\Big) \Big(\expec{\ind{Y \le C} D(\tilde{\ell}_1(\vt_0))\mid X,C}\Big)^T}.
		\end{align*}
		Since $D(\tilde{\ell}_1(\vt_0)) = D(\log(f(Y\mid \vt_0,X)))$ is $\sigma(X,Y)$-measurable and \revised{$C$ is independent of $Y$ given $X$}, we can write
		\begin{align*}
			\expec{\ind{Y \le C} D(\tilde{\ell}_1(\vt_0)) \big(D(\tilde{\ell}_1(\vt_0))\big)^T} &= \expec{\expec{\ind{Y \le C}\mid X,Y} D(\tilde{\ell}_1(\vt_0)) \big(D(\tilde{\ell}_1(\vt_0))\big)^T}\\
			&= \expec{\qty(1-\revised{G(Y- \mid X)}) D(\tilde{\ell}_1(\vt_0)) \big(D(\tilde{\ell}_1(\vt_0))\big)^T}.
		\end{align*}
		Altogether, this yields
		\begin{align*}
			\cov{s_1(\vt_0)} &= \expec{\qty(\revised{\bar{G}(Y- \mid X)}) D(\tilde{\ell}_1(\vt_0)) \big(D(\tilde{\ell}_1(\vt_0))\big)^T}\\
			&\quad + \expec{\qty(\bar{F}(C\mid \vt_0,X))^{-1} \Big(\expec{\ind{Y \le C} D(\tilde{\ell}_1(\vt_0))\mid X,C}\Big) \Big(\expec{\ind{Y \le C} D(\tilde{\ell}_1(\vt_0))\mid X,C}\Big)^T}\\
			&= \Sigma.
		\end{align*}
		
		By the CLT, it follows that $n^{-1/2} s_n(\vt_0)$ weakly converges to a multivariate normal distribution with mean zero and covariance matrix $\Sigma$. Furthermore, $n^{-1} D(s_n(\vt_0))$ converges by the SLLN and equation (\ref{eq:cov_exp}) almost surely to $-\Sigma$. Finally, equation (\ref{eq:lemma_rep}) together with \citet[Lemma 5.54]{dikta2021bootstrap} complete the proof. Note that assumption (\ref{cond:mle_norm__sigma}) guarantees that $\Sigma$ is invertible.
	\end{proof}
	
	\section{Additional simulation results}
	\label{sec:app_sim}
\begin{figure}[H]
		\includegraphics[width=\textwidth]{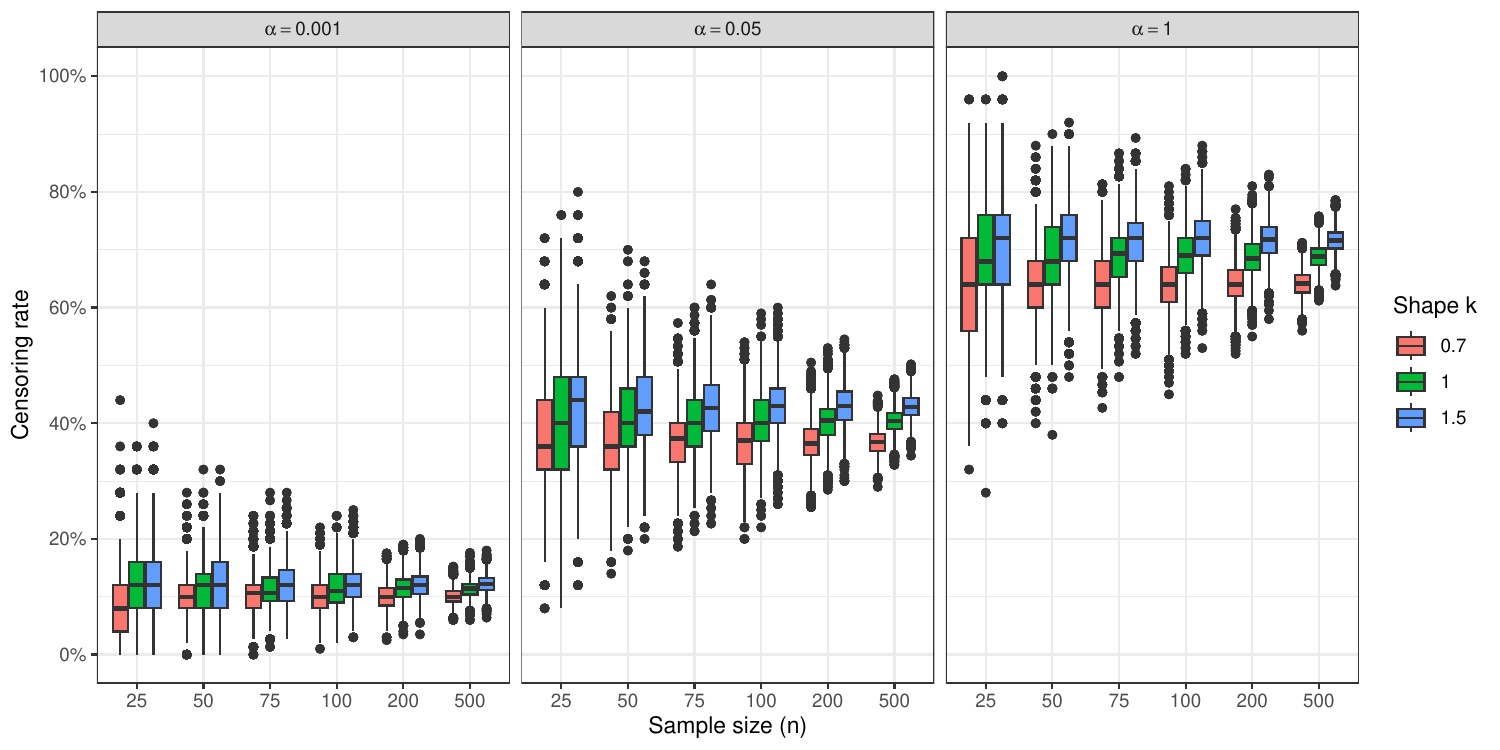}
		\caption{\revised{Distribution of the resulting censoring rates for different choices of $\alpha$, $k$ and $n$. The results are based on the same $500$ simulation replications used for the analysis in the main manuscript.}}
		\label{fig:cens_rates}
	\end{figure}

\begin{figure}[H]
		\centering
		\includegraphics[width=0.8\textwidth]{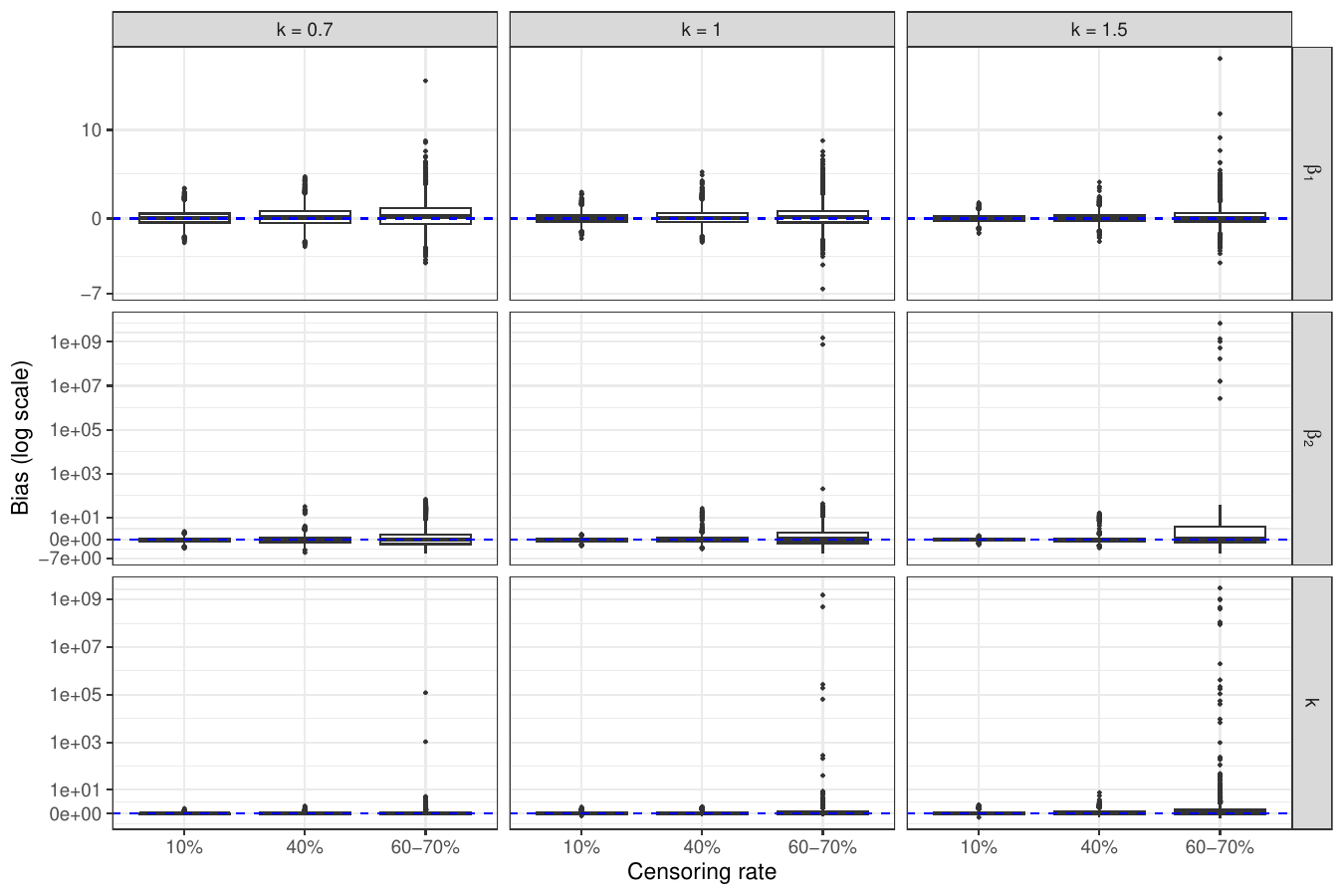}
		\caption{\revised{Empirical results for a very small sample size $n=25$. Same simulation setup as described in section~\ref{sec:sim}. Different censoring times are achieved by varying $\alpha \in 0.001$, $0.05$, $1$.}}
		\label{fig:n25}
	\end{figure}

	\revised{Figure~\ref{fig:cens_rates} shows the distribution of the censoring rates in the simulation study for different choices of $\alpha$. The shape parameter $k$ has only minor influence on the resulting censoring rates, whereas smaller sample sizes lead to greater variability, as expected. 
		
		Figure~\ref{fig:n25} illustrates the empirical bias of the MLE using a sample size of $n=25$ in the simulation study. Even for this very small sample size, the estimates remain essentially unbiased. However, a high censoring rate of $60$--$70\%$ results in a very large variance of the estimates. This is expected, as only about $7$--$10$ event times are observed on average in this case.}
	
	\section{Additional results for the PBC data}
	\label{sec:app_pbc}
	
	Table~\ref{tab:PBC2} reports the estimates and corresponding 95\% confidence intervals for the intercept $\beta_0$ and shape parameter $k$ obtained from fitting a Weibull model to the PBC dataset.

	\begin{center}
		\captionof{table}{Estimates and 95\% confidence intervals for the intercept $\beta_0$ and shape parameter $k$.}
		\label{tab:PBC2}
		\centering
		\begin{tabular}{lcc}
			\toprule
			Method & $\beta_0$ & $k$\\
			\midrule
			Ours & 7.31 (6.00, 8.62) & 1.44 (1.25, 1.62)\\
			Survreg & 7.40 (6.10, 8.71) & 1.44 (1.26, 1.62)\\
			Bootstrap & 7.32 (6.04, 8.65) & 1.46 (1.28, 1.69)\\
			\bottomrule
		\end{tabular}
	\end{center}

	\section{Proofs of lemmas and corollaries}
	
	\begin{proof}[Proof of Lemma \ref{lem:KL1}]
		For part (i), we follow the approach from \citet{stute1992strong} using a Taylor expansion. Let $k(t)=t \log(t)$ for $t>0$ and $k(0)=0$. %
		The Taylor expansion of $k$ at $t=1$ is given by 
		$k(t) = (t-1) + \frac{1}{2\Delta} (t-1)^2$
		with $\Delta$ between $t$ and $1$. Now, let $t_1=\frac{f(y\mid \vt_1, x)}{f(y\mid \vt_2, x)}$ and $t_2 = \frac{1-F(c\mid \vt_1, x)}{1-F(c\mid \vt_2, x)}$. Then,
		\begin{align}
			K_{G,H}(\vt_1, \vt_2) &= \int \int \Bigg[ \bigg[ \int_{-\infty}^c \log \left( \frac{f(y\mid \vt_1, x)}{f(y\mid \vt_2, x)} \right) \frac{f(y\mid \vt_1, x)}{f(y\mid \vt_2, x)} f(y\mid \vt_2,x) \nu(dy) \bigg] \nonumber\\
			&\qquad\quad + \log \left( \frac{1-F(c\mid \vt_1, x)}{1-F(c\mid \vt_2, x)} \right) \frac{1-F(c\mid \vt_1, x)}{1-F(c\mid \vt_2, x)} (1-F(c\mid \vt_2,x)) \Bigg] \revised{G(dc\mid x) H(dx)} \nonumber\\
			&= \int \int \Bigg[ \bigg[ \int_{-\infty}^c k(t_1) f(y\mid \vt_2,x) \nu(dy) \bigg] + k(t_2) (1-F(c\mid \vt_2,x)) \Bigg] \revised{G(dc\mid x) H(dx)} \nonumber\\
			&= \int \int \Bigg[ \bigg[ \int_{-\infty}^c \left( (t_1 - 1) + \frac{1}{2\Delta_1} (t_1-1)^2 \right) f(y\mid \vt_2,x) \nu(dy) \bigg] \nonumber\\
			&\qquad\quad + \left( (t_2 - 1) + \frac{1}{2\Delta_2} (t_2-1)^2 \right) (1-F(c\mid \vt_2,x)) \Bigg] \revised{G(dc\mid x) H(dx)} \label{eq:K_GH}
		\end{align}
		with $\Delta_1$ between $t_1$ and $1$ and $\Delta_2$ between $t_2$ and $1$. 
		The leading terms add up to zero since
		\begin{align*}
			&\int \int \int_{-\infty}^c \left( \frac{f(y\mid \vt_1, x)}{f(y\mid \vt_2, x)} - 1 \right) f(y\mid \vt_2,x) \nu(dy) \revised{G(dc\mid x) H(dx)}\\
			&\quad = \int \int F(c\mid \vt_1, x) - F(c\mid \vt_2,x) \revised{G(dc\mid x) H(dx)}\\
			&\quad = - \int \int \left( \frac{1-F(c\mid \vt_1, x)}{1-F(c\mid \vt_2, x)} - 1 \right) (1-F(c\mid \vt_2,x)) \revised{G(dc\mid x) H(dx)}.
		\end{align*}
		The remainder, on the other hand, is non-negative as
		\begin{align*}
			f(y\mid \vt_1, x) \ge 0,~ f(y\mid \vt_2,x) \ge 0 &\quad\Rightarrow\quad t_1 = \frac{f(y\mid \vt_1, x)}{f(y\mid \vt_2, x)} \ge 0,\\
			1-F(c\mid \vt_1, x) \ge 0,~ 1-F(c\mid \vt_2,x) \ge 0 &\quad\Rightarrow\quad t_2 = \frac{1-F(c\mid \vt_1, x)}{1-F(c\mid \vt_2, x)} \ge 0,
		\end{align*}
		and hence $\Delta_1,~ \Delta_2 \ge 0$. Altogether, we get $K_{G,H}(\vt_1, \vt_2) \ge 0$.
		
		For part (ii), first assume that 
		$\int \int \indset{A_{\vt_1,\vt_2}} (x,y) f(y\mid \vt_1,x) \, \nu(dy) \, H(dx) = 0.$
		For the first part of the remainder in equation \eqref{eq:K_GH}, using Fubini's Theorem, we have
		\begin{align*}
			&\int \int \int_{-\infty}^c \frac{1}{2 \Delta_1} \qty( \frac{f(y\mid \vt_1,x)}{f(y\mid \vt_2,x)} -1 )^2 f(y\mid \vt_2,x) \nu(dy) \revised{G(dc\mid x) H(dx)}\\
			&\quad = \int \int \frac{1}{2 \Delta_1} \qty( \frac{f(y\mid \vt_1,x)}{f(y\mid \vt_2,x)} -1 )^2 f(y\mid \vt_2,x) \qty(1-G(y-\mid x)) \nu(dy) H(dx)\\
			&\quad = \int \int \frac{1}{2 \Delta_1} \qty( \frac{f(y\mid \vt_1,x)}{f(y\mid \vt_2,x)} -1 )^2 \qty(1-G(y-\mid x)) \indset{A_{\vt_1,\vt_2}}(x,y) \frac{f(y\mid \vt_2,x)}{f(y\mid \vt_1,x)} f(y\mid \vt_1,x) \, \nu(dy) \, H(dx)  = 0.
		\end{align*}
		
		Next, for every $x$, we write $\qty{y \mid G(y-\mid x) < 1} = (-\infty, b(x)]$ where $b(x)$ is the smallest point at which $G(\cdot\mid x)$ becomes one, and $(-\infty, \infty)$ if no such points exist. Using Jensen's inequality for the second part of the remainder in equation~\eqref{eq:K_GH} yields
		\begin{align}
			0 &\le \int \int \frac{1}{2\Delta_2} \qty(\frac{1-F(c\mid \vt_1, x)}{1-F(c\mid \vt_2, x)} - 1)^2 \Big(1-F(c\mid \vt_2, x)\Big) \revised{G(dc\mid x) H(dx)} \nonumber\\
			&= \int_{-\infty}^b \int \frac{1}{2\Delta_2} \Big(F(c\mid \vt_2, x) - F(c\mid \vt_1, x)\Big)^2 \qty(\frac{1}{1-F(c\mid \vt_2, x)}) \revised{G(dc\mid x) H(dx)} \nonumber\\
			&\le \int_{-\infty}^b \int \frac{1}{2\Delta_2} \int_{-\infty}^c \Big(f(y\mid \vt_2, x) - f(y\mid \vt_1, x)\Big)^2 \nu(dy) \qty(\frac{1}{1-F(c\mid \vt_2, x)}) \revised{G(dc\mid x) H(dx)} \nonumber\\
			&= \int_{-\infty}^b \int \int \frac{1}{2\Delta_2} \mathbb{I}_{\{y \le c\}} \frac{\Big(f(y\mid \vt_2, x) - f(y\mid \vt_1, x)\Big)^2}{f(y\mid \vt_1,x)} \qty(\frac{1}{1-F(c\mid \vt_2, x)}) \label{eq:first_remainder}\\
			&\qquad\qquad\qquad \cdot \mathbb{I}_{A_{\vt_1,\vt_2}} f(y\mid \vt_1,x) \nu(dy) \revised{G(dc\mid x) H(dx)} \nonumber\\
			& = 0. \nonumber
		\end{align}
		
		For the converse, assume that $\int \int \mathbb{I}_{A_{\vt_1,\vt_2}} (x,y) f(y\mid \vt_1,x) \, \nu(dy) \, H(dx) > 0$. In that case, it can be easily seen from equation \eqref{eq:first_remainder} that the first part of the remainder in equation \eqref{eq:K_GH} will be greater than zero. Together with part (i) of this lemma, it follows that $K_{G,H}(\vt_1, \vt_2) > 0$.
	\end{proof}

	\begin{proof}[Proof of Lemma \ref{lem:KL2}]
		Let $\vt \in \Theta \backslash \{ \vt_0 \}$, then by assumption (ii) and \fref{lem:KL1},
		$K_{G,H}(\vt_0, \vt) > 0$ %
		implying that $L_{G,H}(\vt_0, \vt) < L_{G,H}(\vt_0, \vt_0)$.
	\end{proof}
	
	\begin{proof}[Proof of \fref{cor:asym_norm}]
		$\tilde{\Sigma}$ being positive definite means that for all $z \ne 0$
		\begin{align*}
			0 < z^T \tilde{\Sigma} z = \expec{z^T D(\tilde{\ell}_1(\vt_0))\qty(D(\tilde{\ell}_1(\vt_0)))^T z} = \expec{\dprod{z,D(\tilde{\ell}_1(\vt_0))}^2},
		\end{align*}
		which implies that $\mathbb{P}(\langle z,D(\tilde{\ell}_1(\vt_0))\rangle^2 > 0) > 0$.
		Further, by assumption \eqref{cond:mle_norm__censoring}, we have
		\begin{align*}
			\prob{\revised{\bar{G}(Y-\mid X)} \dprod{z,D(\tilde{\ell}_1(\vt_0))}^2 > 0} &= \prob{\revised{\bar{G}(Y-\mid X)} > 0, \dprod{z,D(\tilde{\ell}_1(\vt_0))}^2 > 0}\\
			&= \prob{G(Y-) < 1, \dprod{z,D(\tilde{\ell}_1(\vt_0))}^2 > 0}\\
			&= \prob{\dprod{z,D(\tilde{\ell}_1(\vt_0))}^2 > 0} > 0,
		\end{align*}
		and hence
		\begin{equation*}
			0 < \expec{\revised{\bar{G}(Y-\mid X)} \dprod{z,D(\tilde{\ell}_1(\vt_0))}^2} = z^T \expec{\revised{\bar{G}(Y-\mid X)} D(\tilde{\ell}_1(\vt_0)) \qty(D(\tilde{\ell}_1(\vt_0)))^T}z.
		\end{equation*}
		As this holds for any $z \ne 0$, it follows that the first summand of $\Sigma$, namely 
		\begin{equation*}
			\Sigma_1 = \expec{\revised{\bar{G}(Y-\mid X)} D(\tilde{\ell}_1(\vt_0)) \qty(D(\tilde{\ell}_1(\vt_0)))^T}
		\end{equation*}
		is positive definite. Since the second summand of $\Sigma$, given by
		\begin{equation*}
			\Sigma_2 = \expec{\qty(\bar{F}(C\mid \vt_0,X))^{-1} \Big(\expec{\ind{Y \le C} D(\tilde{\ell}_1(\vt_0))\mid X,C}\Big) \Big(\expec{\ind{Y \le C} D(\tilde{\ell}_1(\vt_0))\mid X,C}\Big)^T},
		\end{equation*}
		represents the expected value of $AA^T$ for some matrix $A$ and is thus positive semi-definite, the sum $\Sigma = \Sigma_1 + \Sigma_2$ is positive definite as required by assumption (\ref{cond:mle_norm__sigma}) of the previous theorem.
	\end{proof}
	
	\begin{proof}[Proof of \fref{cor:lin_rep}]
		Using the result established in the proof of \fref{thm:asym_norm} that by the SLLN $n^{-1}D(s_n(\vt_0))$ converges a.s.\ to $-\Sigma$ and hence $n^{-1} D(s_n(\vt_0)) + \Sigma = \opI$, we can rewrite equation (\ref{eq:lemma_rep}) as
		\begin{align*}
			-n^{-1/2} s_n(\vt_0) = (-\Sigma + \opI)n^{1/2}(\vtn-\vt_0).
		\end{align*}
		Next, we want to multiply by the inverse of $(-\Sigma + \opI)$ on both sides. For that, note that, according to e.g.\ \citet{miller1981inverse}, the inverse of the sum of two matrices is given by
		\begin{align*}
			(A+B)^{-1} = A^{-1} - (1+tr(BA^{-1}))^{-1} A^{-1} B A^{-1}.
		\end{align*}
		It follows that $(-\Sigma + \opI)^{-1} = -\Sigma^{-1} + \opI$ and hence
		\begin{align*}
			n^{1/2}(\vtn-\vt_0) %
			&= n^{-1/2} \Sigma^{-1} s_n(\vt_0) + \opI.
		\end{align*}	
		The second and third assertion follow directly from $\expec{s_n(\vt_0)}=0$ and $\cov{s_n(\vt_0)} = \Sigma$ as shown in the proof of \fref{thm:asym_norm}.
	\end{proof}
\end{document}